\documentclass[reqno,centertags, 11pt]{amsart}
\usepackage{comment}
\usepackage{color}
\usepackage{amsmath,amsfonts,amssymb,amscd,amsthm,amsbsy, epsf,calc,enumerate,cite}
\usepackage{epsfig,epstopdf,esint,psfrag}
\usepackage{graphics,graphicx,gensymb}
\usepackage{amssymb,amsmath,amsfonts}
\advance\hoffset by 0.4mm
\usepackage{latexsym}

\def\@abssec#1{\vspace{.05in}\footnotesize \parindent .2in
{\bf #1. }\ignorespaces}
\newtheorem{theorem}{Theorem}[section]

\newtheorem{lemma}[theorem]{Lemma}

\def \Rm {\mathbb R}
\def \Tm {\mathbb T}

\newcommand{\eps}{\varepsilon}

\newcommand{\be}{\mathbf e}

\newcommand{\lb}{\label}

\allowdisplaybreaks \numberwithin{equation}{section}

\renewcommand{\be}{\begin{equation}}
\newcommand{\ee}{\end{equation}}

\newcommand{\om}{\omega}
\newcommand{\OO}{\mathcal O}

\title[Singularity formation in fluids]{Small scales and singularity formation in fluid dynamics}
\author{Alexander Kiselev}
\thanks{Department of
Mathematics, Duke University, Durham NC 27708, USA;
email: kiselev@math.duke.edu}

\begin{document}


\begin{abstract}
We review recent advances in understanding singularity and small scales formation in solutions of fluid dynamics equations.
The focus is on the Euler and surface quasi-geostrophic (SQG) equations and associated models.
\end{abstract}


\maketitle

\section{Introduction}\label{intro}

Fluids are all around us, and attempts to mathematically understand fluid motion go back many centuries. Anyone witnessing a dramatic phenomena like
tornado or hurricane or even an everyday river flow or ocean wave breaking can easily imagine the complexity of the task. There has been
tremendous accumulation of knowledge in the field, yet it is remarkable that some of the fundamental properties of key equations of fluid mechanics
remain poorly understood.

A special role in fluid mechanics is played by the incompressible Euler equation, first formulated in 1755 \cite{Euler}. Amazingly,
it appears to be the second partial differential equation ever derived (the first one is wave equation derived by D'Alembert 8 years earlier).
The incompressible Euler equation describes motion of an inviscid, volume preserving fluid; fluid with such properties is often called ``ideal".
The Euler equation is a nonlinear and nonlocal system of PDE, with dynamics near a given point depending on the flow field over the entire region filled with fluid.
This makes analysis of these equations exceedingly challenging, and the array of mathematical methods applied to their study has been very broad.

The basic purpose of an evolution PDE is solution of Cauchy problem: given initial data, find a solution that can then be used for prediction
of the modelled system. This is exactly how weather forecasting works, or how new car and airplane shapes are designed.
Therefore, one of the first questions one can ask about a PDE is existence and uniqueness of solutions in appropriate
functional spaces. The PDE is called globally regular if there exists a unique, sufficiently smooth solution for reasonable classes of initial data.
On the other hand, singularity formation - meaning that some quantities associated with fluid motion become infinite - can indicate spontaneous
creation of intense fluid motion. Singularities are also important to understand since they may indicate
potential breakdown of the model, may lead to loss of uniqueness and predictive power, and are very hard to resolve computationally.
More generally, one can ask a related and broader question about creation of small scales in fluids - coherent structures that vary sharply in space and time,
and contribute to phenomena such as turbulence (see e.g. \cite{Eyink} for further references).

For the Euler equation, the global regularity vs finite time blow up story is very different depending on the dimension. Let us recall that the
incompressible Euler equation in a domain $D \subset \Rm^d,$ $d=2$ or $3$ with natural no boundary penetration conditions is given by
\begin{equation}\label{euler}
\left. \partial_t u + (u \cdot \nabla) u = \nabla p, \,\,\,\,\,\,\nabla \cdot u =0,\,\,\,\,\,\,u \cdot n \right|_{\partial D}=0,
\end{equation}
along with the initial data $u(x,0)=u_0.$ Here $u(x,t)$ is the vector field describing fluid velocity, $p$ is pressure, and $n$ is the normal at the boundary $\partial D.$ The equation \eqref{euler} is just the second Newton's law written for ideal fluid.
The difference between dimensions becomes clear if we rewrite the
equation in vorticity $\omega = {\rm curl} u:$
\begin{equation}\label{eulervort}
\partial_t \omega +(u \cdot \nabla) \omega = (\omega \cdot \nabla) u,\,\,\,\,\,\,\omega(x,0)=\omega_0(x),
\end{equation}
along with Biot-Savart law which allows to recover $u$ from $\omega.$ For instance, in the case of a smooth domain $D \subset \Rm^2$ one gets $u = \nabla^\perp (-\Delta_D)^{-1}\omega$ -
where $\nabla^\perp = (\partial_{x_2}, -\partial_{x_1})$ and $\Delta_D$ is the Dirichlet Laplacian.
In the form \eqref{eulervort}, one can observe that the term $(\omega \cdot \nabla) u$ on the right hand side vanishes in dimension two.
The resulting equation conserves the
$L^\infty$ norm and in fact any $L^p$ norm of a regular solution, which helps prove global regularity. This result has been known since 1930s works
by Wolibner \cite{Wolibner} and Holder \cite{Holder}. We will focus on the 2D Euler equation in Section~\ref{doubleexp} below.
We note that another feature of the Euler equation made obvious by the vorticity representation \eqref{eulervort} (in any dimension) is nonlocality:
the inverse Laplacian in the Biot-Savart law is a manifestly nonlocal operator, involving integration over the entire domain.

In three dimensions, the ``vortex-stretching" term $(\omega \cdot \nabla) u$ is present and can affect the intensity of vorticity.
Local well-posedness results in a range of natural spaces are well known; one can consult \cite{MB} or \cite{MP} for proofs and further references.
However, the global regularity vs finite time singularity formation question remains open. In fact, this problem for the Euler equation is a close
relative of the celebrated Clay Institute Millennium problem on the 3D Navier-Stokes equation \cite{Feff}. Indeed, the Navier-Stokes equation only differs from
\eqref{euler} by the presence of a regularizing, linear term $\Delta u$ on the right hand side modeling viscosity (and by different boundary conditions
if boundaries are present). The nonlinearity - the principal engine of possible singular growth - is identical in both equations.

A host of numerical experiments sought to discover a scenario for singularity formation in solutions of 3D Euler equation (see e.g. \cite{BNW,BorPel,EShu,GraSid,Hou1,Kerr,LPTW,OG,PG,PumSig},
or a detailed review by Gibbon \cite{Gibbon} where more references can be found). In the analytic direction, a complete review would be too broad to
attempt here. Let us mention the classical work of Beale, Kato and Majda \cite{BKM} on criteria for global regularity, papers \cite{Caf,Caf1} on singularities for
complex-valued solutions to Euler equations, as well as regularity criteria by Constantin, Fefferman and Majda \cite{CFM} and by Hou and collaborators \cite{Hou1,Hou11,DHY} which involve
more subtle geometric conditions sufficient for regularity. See Peter Constantin \cite{PC12} for more history and analytical aspects of this problem.

There have been several recent developments in classical problems on regularity and solution estimates for the fundamental equations of fluid mechanics.
First, Hou and Luo produced a new set of careful numerical experiments suggesting finite time singularity formation for solutions of 3D Euler equation \cite{HouLuo}.
The scenario of Hou and Luo is axi-symmetric. Very fast vorticity growth is observed at a ring of hyperbolic stagnation points of the flow
located on the boundary of a cylinder. None of the available regularity criteria such as \cite{CFM,Hou1,DHY,Hou11} seem to apply to scenario's geometry.
The scenario has a close analog for the 2D inviscid Boussinesq system, for which the question of global regularity is also open; it is listed as one of the
``eleven great problems of hydrodynamics" by Yudovich \cite{Yud2000}. More details on the scenario will be provided in Section~\ref{hyp}.

The hyperbolic stagnation point on the boundary scenario also leads to interesting phenomena in solutions of the 2D Euler equation. The best known upper bound on the growth of
the gradient of vorticity, as well as higher order Sobolev norms, is double exponential
in time:
\begin{equation}\label{ub1} \|\nabla \omega(\cdot, t)\|_{L^\infty} \leq \left(1+\|\nabla \omega_0\|_{L^\infty}\right)^{\exp(C\|\omega_0\|_{L^\infty}t)}. \end{equation}
This result has appeared explicitly in \cite{YudDE}, though related bounds can be traced back to \cite{Wolibner}.
The question whether such upper bounds are sharp has been open for a long time.
In a joint work with Sverak \cite{KS}, we provided an example of smooth initial data in the disk such that the corresponding solution 
exhibits double exponential growth in the gradient of vorticity for all times, establishing qualitative sharpness of \eqref{ub1}.
The construction is based on the hyperbolic point at the boundary scenario, and will be described in more detail in Section~\ref{doubleexp}.

Further attempts to rigorously understand the Hou-Luo scenario involved construction of 1D and 2D models retaining some of the analytical structure
of the original problem. We will discuss some of these models in Section~\ref{models}.


The surface quasi-geostrophic (SQG) equation is similar to the 2D Euler equation in vorticity form, but is more singular:
\begin{equation}\label{sqg} \partial_t \omega + (u \cdot \nabla) \omega =0, \,\,\, u = \nabla^\perp (-\Delta)^{-1+\alpha}\omega, \,\,\,\alpha=1/2,\,\,\,\omega(x,0)=\omega_0(x). \end{equation}
The value $\alpha=0$ corresponds to the 2D Euler equation, while $0 < \alpha < \frac12$ is called the modified SQG range. The SQG and modified SQG equations come from atmospheric science.
They model evolution of temperature near the surface of a planet and can be derived by formal asymptotic analysis from a larger system of rotating 3D Navier-Stokes equations
coupled with temperature equation through buoyancy force \cite{Held,Majda,Ped,PHS}. In mathematical literature, the SQG equation was first considered by Constantin, Majda and Tabak
\cite{CMT}, where a parallel between the structure of the SQG equation and the 3D Euler equation was drawn. The SQG and modified SQG equations are perhaps simplest looking equations
of fluid mechanics for which the question of global regularity vs finite time blow up remains open. The equation \eqref{sqg} can be considered with smooth initial data, but another important
class of initial data is patches, where $\theta_0(x)$ equals linear combination of characteristic functions of some disjoint domains $\Omega_j(0).$ The resulting evolution yields time
dependent regions $\Omega_j(t)$.
The regularity question in this context addresses the regularity class of the boundary $\partial \Omega_j(t)$ and lack of contact between different components.
Existence and uniqueness of patch solution for 2D Euler equation
follows from Yudovich theory \cite{Yudth,MB,MP}.
The global regularity question has been settled affirmatively by Chemin \cite{c} 
(Bertozzi and Constantin \cite{bc} provided a different proof). For the SQG and modified SQG equations patch dynamics is harder to set up.
Local well-posedness has been shown by Rodrigo in $C^\infty$ class \cite{Rodrigo} and by Gancedo in Sobolev spaces \cite{g} in the whole plane setting.
Numerical simulations by Cordoba, Fontelos, Mancho and Rodrigo \cite{CFMR} suggest that finite time singularities -
in particular formation of corners and different components touching each other - is possible, but rigorous understanding of this phenomena remained missing.
Jointly with Ryzhik, Yao and Zlatos \cite{KRYZ}, \cite{KYZ}, we considered modified SQG and 2D Euler patches in half-plane, with natural no penetration boundary conditions. The initial
patches are regular and do not touch each other but may touch the boundary. We proved a kind of phase transition in this setting: the 2D Euler patches stay globally regular, while for a
range of small $\alpha>0$ some initial data lead to blow up in finite time.
The blow up scenario again involves a stagnation hyperbolic point of the flow on the boundary.
This result will be described in more detail in Section~\ref{patch}.


{\bf Acknowledgement.} This work has been supported by the NSF-DMS award 1712294.

\section{The 3D Euler equation and the 2D Boussinesq system: the hyperbolic scenario}\label{hyp}

\begin{figure}
\begin{center}\label{fig11}
\includegraphics[width=100mm]{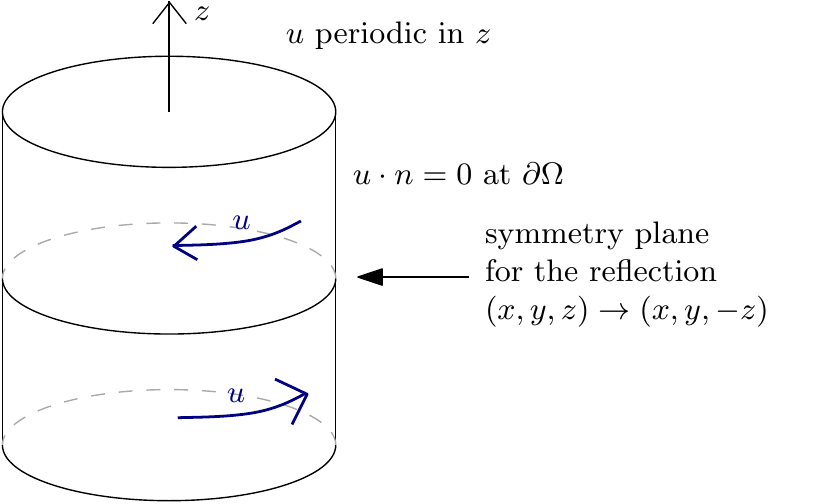}
\caption{The initial data for Hou-Luo scenario}
\end{center}
\end{figure}

In~\cite{HouLuo} the authors study 3D axi-symmetric solutions of incompressible Euler equation with roughly the initial configuration shown on Figure~\ref{fig11}:
only swirl $u^\phi$ is initially non-zero, and it is odd and periodic in $z$ variable.

One of the standard forms of the axi-symmetric Euler equations in the usual cylindrical coordinates $(r,\phi, z)$ is
\begin{subequations}\label{axi-sym}
\begin{align}
\partial_t \left(\frac{\omega^\phi}{r}\right)+
u^r\partial_r \left(\frac{\omega^\phi}{r}\right) + u^z \partial_z\left(\frac{\omega^\phi}{r}\right) & =\partial_z\left(\frac{(ru^\phi)^2}{r^4}\right)\\
\partial_t (ru^\phi)+u^r \partial_r (ru^\phi)+u^z \partial_z (ru^\phi) & = 0\,,
\end{align}
\end{subequations}
with the understanding that $u^r, u^z$ are given from $\om^\phi$ via the Biot-Savart law 
which in the setting of Hou-Luo scenario takes form
\[ u^r = -\frac{\partial_z \psi}{r}, \,\,\, u^z = \frac{\partial_r \psi}{r}, \,\,\, L\psi = \frac{\omega^\phi}{r}, \,\,\,
L\psi = -\frac1r \partial_r \left( \frac1r \partial_r \psi \right) - \frac{1}{r^2} \partial_{zz}^2 \psi. \]

\begin{figure}\label{secflowsfig}
\begin{center}
\includegraphics[width=100mm]{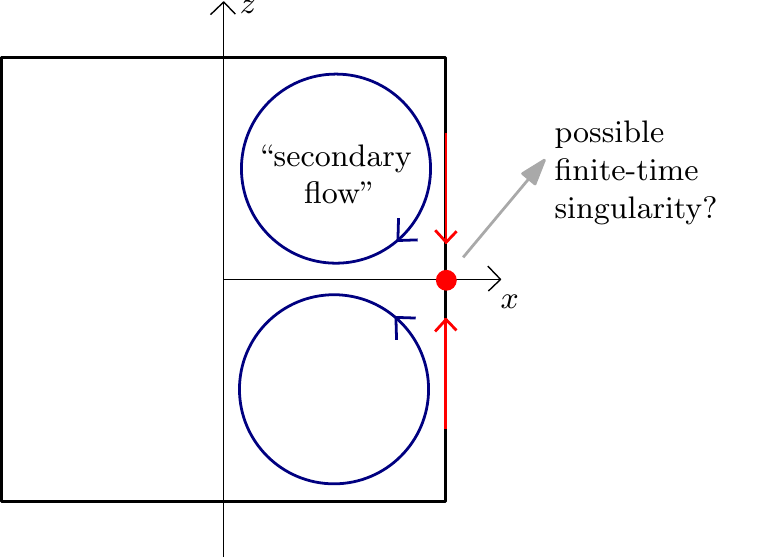}
\caption{The secondary flows in fixed $\phi$ section}
\end{center}
\end{figure}

From \eqref{axi-sym}, it is clear that the swirl will spontaneously generate toroidal rolls corresponding to non-zero $\omega^\phi.$
These are the so-called ``secondary flows",~\cite{prandtl}; its effect on river flows was studied by Einstein \cite{einstein}.
Thus the initial condition leads to the (schematic)
picture in the $xz$--plane shown on Figure~\ref{secflowsfig}, in which we also indicate the point where a conceivable
finite-time singularity (or at least an extremely strong growth of vorticity) is observed numerically.
In the three-dimensional picture, the points with very fast growth form a ring on the boundary of the cylinder.

A somewhat similar scenario can be considered for the 2D inviscid Boussinesq system in a half-space $\Rm^+=\{(x,y)\in\Rm\times (0,\infty)\}$ (or in a flat half-cylinder ${\bf S}^1\times (0,\infty)$),
which we will write in the vorticity form:
\begin{subequations}\label{boussinesq}
\begin{align}
\partial_t \om+u_1 \partial_x \om +u_2 \partial_y \om & = \partial_x \theta\\
\partial_t \theta+u_1 \partial_x \theta + u_2 \partial_y \theta & = 0\,.
\end{align}
\end{subequations}
 Here $u=(u_1,u_2)$ is obtained from $\omega$ by the usual Biot-Savart law $u = \nabla^\perp (-\Delta)^{-1}\omega,$ with appropriate boundary
 conditions on $\Delta,$
and $\theta$ represents the fluid temperature or density.

\begin{figure}\label{Bousfig}
\begin{center}
\includegraphics[width=90mm]{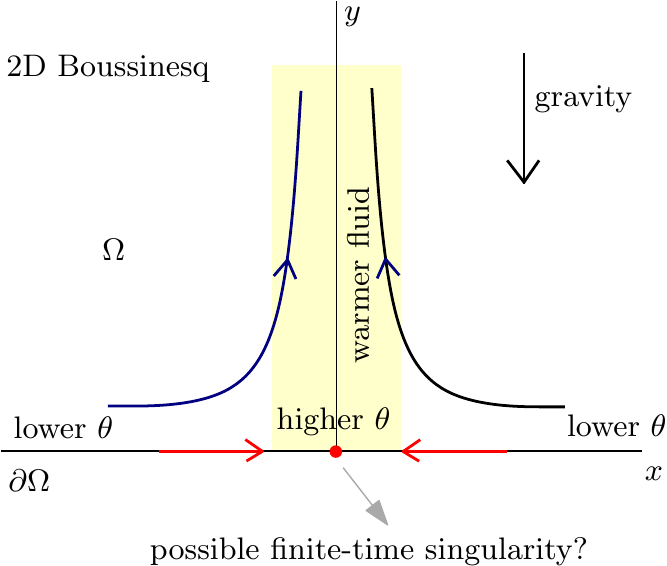}
\caption{The 2D Boussineq singularity scenario}
\end{center}
\end{figure}

It is well-known (see e.g. \cite{MB}) that this system has properties similar to the 3D axi-symmetric Euler \eqref{axi-sym},
at least away from the rotation axis. Indeed, comparing \eqref{axi-sym} with \eqref{boussinesq}, we see that
$\theta$ essentially plays the role of the square of the swirl
component $r u^\phi$ of the velocity field $u$, and $\omega$ replaces $\omega^\phi/r.$
The real difference between the two systems only emerges near the axis of rotation, where the factors of $r$ can conceivably
change the nature of dynamics.
For the purpose of comparison with the axi-symmetric flow, the last picture should be rotated
by $\pi/2$, after which it resembles the picture relevant for~\eqref{boussinesq}, see Figure~\ref{Bousfig}.

In both the 3D axi-symmetric Euler case and in the 2D Boussinesq system case the best chance for possible singularity formation
seems to be at the points of symmetry at the boundary,
which numerical simulations suggest are fixed hyperbolic points of the flow.

\section{The 2D Euler equation}\label{doubleexp}

A reasonable first step to understand the hyperbolic stagnation point on the boundary blow up scenario is to consider the case of constant
density in the Boussinesq system first. Of course, this reduces the system to the 2D Euler equation:
\begin{equation}\label{2devort}
\partial_t \omega + (u \cdot \nabla) \omega =0, \,\,\,u = \nabla^\perp(-\Delta_D)^{-1}\omega, \,\,\,\omega(x,0)= \omega_0(x).
\end{equation}
Here $\Delta_D$ stands for Dirichlet Laplacian; such choice of the boundary condition corresponds to no-penetration property $\left. u \cdot n \right|_{\partial D} =0.$
Of course, for the 2D Euler equation solutions are globally regular.
Let us state this result, going back to 1930s \cite{Wolibner,Holder}.
\begin{theorem}\label{upperb}
Let $D \subset \Rm^2$ be a compact domain with smooth boundary, and $\omega_0(x) \in C^1(D)$
Then there exists a unique smooth solution $\omega(x,t)$ of the equation \eqref{2devort} corrsponding to the initial data $\omega_0$, which moreover satisfies
\begin{equation}\label{upperdexp}
1+\log\left(1+\frac{\|\nabla \omega(x,t)\|_{L^\infty}}{\|\omega_0\|_{L^\infty}}\right) \leq
\left(
1+\log\left(
1+
{
\|\nabla\omega_0\|_{L^\infty}
\over
\|\omega_0\|_{L^\infty}
}
\right)
\right)
\exp(
C\|\omega_0\|_{L^\infty}t
)
\end{equation}
for some constant $C$ which may depend only on the domain $D.$
\end{theorem}

A key ingredient of the proof is the Kato inequality \cite{Kato}: for every $1>\alpha>0,$ we have
\begin{equation}\label{katoeq}
\|\nabla u(x,t)\|_{L^\infty} \leq C(\alpha,D) \|\omega_0\|_{L^\infty}\left(1+\log \frac{\|\omega(x,t)\|_{C^\alpha}}{\|\omega_0\|_{L^\infty}} \right).
\end{equation}
Note that the matrix $\nabla u$ consists of double Riesz transforms $\partial^2_{ij}(-\Delta_D)^{-1} \omega.$ Riesz transforms
are well known to be bounded on $L^p,$ $1< p < \infty$ (see e.g. \cite{Stein}), but the bound fails at the endpoints and we have
to pay a logarithm of the higher order norm to obtain a correct bound. It is exactly the extra $\log$ in \eqref{katoeq} that leads
to the double exponential upper bound as opposed to the single one (see e.g. \cite{KS} for more details).

The question of whether such upper bounds are sharp has been open for a long time. Yudovich \cite{Jud1,Yud2} provided an example showing
 infinite
growth of the vorticity gradient at the boundary of the domain, by constructing an appropriate Lyapunov-type functional. These results were
 further
improved and generalized in \cite{MSY}, leading to description of a broad class of flows with infinite growth in their vorticity gradient.
Nadirashvili \cite{Nad1} proved a more quantitative linear in time lower bound for a ``winding" flow in an annulus.
A variant of the example due to Bahouri and Chemin
\cite{BC} provides singular stationary solution of the 2D Euler equation defined on $\Tm^2$ with fluid velocity which is
just log-Lipschitz in spacial
variables. Namely, if we set $\Tm^2 = [-\pi,\pi) \times [-\pi,\pi),$ the solution is equal to $-1$ in the first quadrant
$[0,\pi) \times [0,\pi)$
and is odd with respect to both coordinate axes. Note that the solution is just $L^\infty$ but existence and uniqueness of solutions in
this class
is provided essentially by Yudovich theory \cite{Yudth}. The origin is a fixed hyperbolic point of the fluid velocity, with $x_1$ being the contracting direction, and
the velocity satisfies $u_1(x_1,0) = \frac{4}{\pi} x_1 \log x_1 +O(x_1)$ for small $x_1.$
The trajectory starting at a point $(x_1,0)$ on a horizontal separatrix will therefore converge to the origin at a double exponential rate
in time. If a smooth passive scalar $\psi$ is advected by a flow generated by singular cross, $\partial_{x_1}\psi$ will grow at double
exponential
rate in time if $\psi(x_1,0)$ is not constant for small $x_1.$ Of course, derivative growth does not make sense for  the singular cross
solution
itself since it is stationary and already singular. But this solution shows a blueprint of how double exponential growth can be conceivably
generated in smooth solution: it needs to approach the discontinuous configuration similar to the singular cross, while at the
same time the solution should be nonzero on the contracting direction. This turns out to be hard to implement, especially without boundary.

In recent years, there has been a series of works by Denisov on this problem. In \cite{Den1}, he constructed an example with superlinear
growth in vorticity gradient of the solution in the periodic case. In \cite{Den2}, he showed that for any time $T$, one can arrange smooth initial
data so that the corresponding
solution will experience double exponential burst of growth over $[0,T].$ The example is based on smoothing out Bahouri-Chemin example,
abandoning odd symmetry to put a ripple on a separatrix, and controlling the resulting solution over finite time interval.

In \cite{KS}, we proved
\begin{theorem}\label{mainthm}
Consider two-dimensional Euler equation on a unit disk $D.$ There exists a smooth initial data $\omega_0$ with
$\|\nabla \omega_0\|_{L^\infty}/\|\omega_0\|_{L^\infty} >1$ such that the corresponding
solution $\omega(x,t)$ satisfies
\begin{equation}\label{main1}
\frac{\|\nabla \omega(x,t)\|_{L^\infty}}{\|\omega_0\|_{L^\infty}}
\geq \left( \frac{\|\nabla \omega_0\|_{L^\infty}}{\|\omega_0\|_{L^\infty}}\right)^{c\exp(c\|\omega_0\|_{L^\infty}t)}
\end{equation}
for some $c>0$ and for all $t \geq 0.$
\end{theorem}
The theorem shows that double exponential growth in the gradient of vorticity can actually happen, so double exponential
upper bound is sharp. As in Hou-Luo blow up scenario, growth happens at a hyperbolic fixed point of the low at the boundary.
The result has been generalized to the case of any compact sufficiently regular domain with
symmetry axis by Xu \cite{Xu}. The question of whether double exponential growth can happen in the bulk of the fluid remains open;
Zlatos \cite{Z3} has improved the techniques behind Theorem~\ref{mainthm} to construct examples of smooth solutions
with exponential growth of $\|\nabla^2\omega\|_{L^\infty}$ in periodic setting. The question of whether double exponential growth
is at all possible in the bulk of the fluid remains wide open.

A key step in the proof is understanding the structure of fluid velocity near the hyperbolic point.
Let \[ D^+ = \{ x \in D \left| x_1 \geq 0. \right. \} \]
We will choose the initial data that is odd in $x_1,$ and $-1 \leq \omega_0(x)< 0$ for $x \in D^+.$
Let us set the origin of our coordinate system at the bottom of the disc, where interesting things will be happening.
Given the symmetry of $\omega,$ we have
\begin{equation}\label{BioS23}
u(x,t) = -\nabla^\perp \int_D G_D(x,y)\omega(y,t)\,dy =
-\frac{1}{2\pi} \nabla^\perp \int_{D^+}
\log\left(\frac{|x-y||\tilde{x}-\bar{y}|}{|x-\bar{y}||\tilde{x}-y|}\right)
\omega(y,t)\,dy,
\end{equation}
where $G_D$ is the Green's function of Dirichlet Laplacian and $\tilde{x} = (-x_1,x_2).$
For each point $(x_1,x_2)\in D^+$, let us introduce the region
\[
Q(x_1,x_2)=\{(y_1,y_2)\in D^+:~~x_1 \leq y_1,~~x_2 \leq y_2\},
\]
and set
\be\lb{omegamain}
\Omega(x_1,x_2,t) = -\frac{4}{\pi}\int_{Q(x_1,x_2)}
\frac{y_1y_2}{|y|^4} \omega(y,t)\,dy_1dy_2.
\ee
Finally, for any $0<\gamma<\pi/2,$ let $\phi$ be the usual polar angle coordinate of point $x$,
and denote
\[ D_1^\gamma = \{ x \in D^+ \left| 0\le\phi\le\pi/2-\gamma \right. \}, \,\,\,
 D_2^\gamma = \{ x \in D^+ \left| \gamma\le\phi\le\pi/2 \right. \}. \]
The following Lemma is crucial for the proof of Theorem~\ref{mainthm}.
\begin{lemma}\label{mainlemma}
Suppose that $\omega_0$ is odd with respect to $x_1.$ Fix a small $\gamma>0.$
There exists $\delta >0$ so that for all
$x \in D^\gamma_1$ such that $|x| \leq \delta$ we have
\begin{equation}\label{velest1}
u_1(x_1,x_2,t) = - x_1 \Omega(x_1,x_2,t) + x_1 B_1(x_1,x_2,t),
\end{equation}
where $|B_1(x_1,x_2,t)| \leq C(\gamma)\|\omega_0\|_{L^\infty}.$

Similarly, for all $x \in D^\gamma_2$ such
that $|x| \leq \delta$ we have
\begin{equation}\label{velest2}
u_2(x_1,x_2,t) = x_2\Omega(x_1,x_2,t) + x_2 B_2(x_1,x_2,t),
\end{equation}
where $|B_2(x_1,x_2,t)| \leq C(\gamma)\|\omega_0\|_{L^\infty}.$
\end{lemma}

Proof of the lemma is based on analysis of \eqref{BioS23}, full details can be found in \cite{KS}.
The terms on the right hand sides of \eqref{velest1}, \eqref{velest2} involving $\Omega$ can be thought of as main
terms in certain regimes. Indeed, observe that as support of the set where $\omega(x,t)\geq c>0$ approaches the origin,
the size of $\Omega$ given by \eqref{omegamain} may grow as a logarithm of this distance. Thus Lemma~\ref{mainlemma}
provides a sort of quantitative version of Bahouri-Chemin log-Lipschitz singular flow for smooth setting.
In the regime where the $\Omega$ terms dominate, the flow trajectories described by \eqref{velest1}, \eqref{velest2} are
close to precise hyperbolas. Another interesting feature of the formulas \eqref{velest1}, \eqref{velest2} is a hidden
comparison principle: the influence region $Q(x_1,x_2)$ tends to be larger for points closer to the origin.
The comparison principle is not precise, but it turns out to be true up to Lipschitz errors. This feature is key
in the construction of the example with double exponential growth.

Let us now sketch the construction of such example. Fix small $\gamma>0,$ and take $\epsilon>0$ smaller than the corresponding
value of $\delta$ from Lemma~\ref{mainlemma}. We also take $\epsilon$ sufficiently small so that $\log \epsilon^{-1}$ is much
larger than the constant $C(\gamma)$ from Lemma~\ref{mainlemma}. 
Take the initial data $\omega_0$ such that $\omega_0(x_1,x_2) =-1$ if $x_1 \geq \epsilon^{10},$
odd with respect to $x_1,$ and satisfying $-1 \leq \omega_0(x) \leq 0$ for $x \in D^+.$ This leaves certain ambiguity in how
we define $\omega_0$ in the region close to $x_2$ axis. As we will see, it does not matter for the construction exactly how $\omega_0$ is defined
there.

The first observation is that with such choice of the initial data,
\be\lb{Oest11}
\Omega(x,t) \geq C\log \epsilon^{-1}\ee for every $x \in D^+$
with $|x| \leq \epsilon.$ Indeed, due to incompressibility of the flow the measure of the points $x \in D^+$ where $\omega(x,t) >-1$
does not exceed $2\epsilon$ for any time $t.$ It is then straightforward to show that even if all these points
are pushed by the dynamics into the region where the size of the kernel in \eqref{omegamain} is maximal, the estimate \eqref{Oest11}
still holds.

Next, for $0<x_1'<x_1''<1$ we  denote
\be\lb{xx1}
\OO(x_1',x_1'')=\left\{
(x_1,x_2)\in D^+\,,\,\,x_1' \leq x_1 \leq x_1''\,,\,\,
x_2 \leq x_1\right\}\,.
\ee
We would like to analyze the evolution in time of the region $\OO_{\epsilon^{10}\epsilon}.$
For this purpose, for $0<x_1<1$ we let
\be\lb{xx2}
{\underline u}_1(x_1,t) \quad =\quad  \min_{(x_1,x_2)\in D^+\,,\,x_2 \leq x_1} u_1(x_1,x_2,t)
\ee
and
\be\lb{xx3}
\,\,{\overline u}_1(x_1,t) \quad = \quad \max_{(x_1,x_2)\in D^+\,,\, x_2 \leq x_1} u_1(x_1,x_2,t)\,.
\ee
Define
$a(t)$ by
\be\lb{xx4}
a'= \overline u_1(a,t)\,,\quad a(0)=\epsilon^{10}\,
\ee
and  $b(t)$  by
\be\lb{xx5}
b' = \underline u_1(b,t)\,,\quad b(0)=\epsilon.
\ee
Let
\be\lb{xx6}
\OO_t=\OO(a(t),b(t))\,.
\ee

We claim that $\omega(x,t)=-1$ for every $x \in \OO_t,$ and every $t \geq 0.$ Indeed, if this is not the case, then there must exist
some time $s$ and a trajectory $\Phi_s(y)$ such that $y \notin \OO_0,$ and $\Phi_s(y) \in \partial \OO_s$
for the first time. But the trajectory $\Phi_s(y)$ cannot enter $\OO_s$ through the boundary of $D$ due to the no-penetration boundary condition.
It is also not hard to see it cannot enter at the $x_1=a(s)$ or $x_1=b(s)$ pieces of $\partial \OO_s$ due to the definition of
$a(t),b(t)$ and $\overline u,$ $\underline u.$ This leaves the diagonal $x_1=x_2.$ However, due to Lemma~\ref{mainlemma} and the estimate
\eqref{Oest11}, we have that
\be\lb{diag1}
\frac{C\log \epsilon^{-1} - C(\gamma)}{C\log \epsilon^{-1} + C(\gamma)} \leq \frac{-u_1(x_1,x_1,t)}{u_2(x_1,x_1,t)} \leq \frac{C\log \epsilon^{-1} + C(\gamma)}
{C\log \epsilon^{-1} - C(\gamma)}
\ee
for every $t \geq 0$ and $x_1 \leq \epsilon.$ Due to choice of $\epsilon,$ we have that the ratio $-u_1/u_2$ on the diagonal part of the boundary of $\OO_s$
is close to $1$. Thus the vector field $u$ points outside of the region $\OO_s$ on the diagonal part of the boundary
at all times and so the trajectory cannot enter through the diagonal either.

Next, we are going to estimate how quickly $a(t)$ approaches the origin. The constant $C$ below depends only on $\gamma$ and may change from line to line.
By Lemma~\ref{mainlemma}, we have
\[ \underline u_1(b(t),t) \geq  -b(t)\, \Omega(b(t), x_2(t))-Cb(t), \]
for some $x_2(t)\le b(t)\,,\,\,(x_2(t),b(t))\in  D^+$
as  $\|\omega(x,t)\|_{L^\infty} \leq 1$ by our choice of the initial datum $\om_0$.
A straightforward calculation shows that \[ \Omega(b(t), x_2(t)) \leq \Omega(b(t), b(t)) + C. \]
Thus we get
\begin{equation}\label{keyestCD}
\underline u_1(b(t),t) \geq -b(t)\, \Omega(b(t),b(t))-Cb(t).
\end{equation}
Similarly,
\[ \overline u_1(a(t),t) \leq -a(t) \,\Omega(a(t), \tilde x_2 (t)) \leq -a(t) \,\Omega(a(t), 0)+Ca(t). \]
These estimates establish a form of comparison principle, up to Lipschitz errors, of the fluid velocities of the front and
back of the region $\OO_{a(t),b(t)}.$

Now observe that
\be\lb{omeulest} \Omega(a(t), 0) \geq -\frac{4}{\pi} \int_{\OO_t} \frac{y_1y_2}{|y|^4} \omega(y,t)\,dy_1dy_2 + \Omega(b(t), b(t)). \ee
Since $\omega(y,t)=-1$ on $\OO_t,$ a direct estimate shows that the integral in \eqref{omeulest}
is bounded from below by  $\kappa(-\log a(t) + \log b(t) )$  for some $\kappa>0.$ 
Applying the above estimates to evolution of $a(t)$ and $b(t)$ we obtain
\begin{equation}\label{final1}
\frac{d}{dt} \left( \log a(t) - \log b(t) \right) \leq \kappa \left(\log a(t)-\log b(t)\right) +C.
\end{equation}
Applying Gronwall lemma and choosing $\epsilon$ small enough leads to $\log a(t) \leq \exp(\kappa t) \log \epsilon.$
To arrive at \eqref{main1}, it remains to note that
we can arrange $\|\nabla \omega_0\|_{L^\infty} \lesssim \epsilon^{-10}.$


\section{The one-dimensional models}\label{models}

One-dimensional models in fluid mechanics have a long history. We briefly review some of the results most relevant to our narrative.
In the context of modeling finite time blow and global regularity issues,
Constantin, Lax and Majda \cite{CLM} considered the model
\be\label{clm}
\partial_t \omega = \omega H \omega, \,\,\,\omega(x,0) = \omega_0(x),
\ee
where $H$ is the Hilbert transform, $H \omega(x,t) = \frac{1}{\pi} P.V. \int_{\Rm} \frac{\omega(y,t)}{x-y}\,dy.$
The equation \eqref{clm} is designed to model the vortex stretching term on the right hand side of \eqref{eulervort}; the advection
term is omitted. Surprisingly, the model \eqref{clm} is explicitly solvable due to special properties of Hilbert transform.
Finite time blow up happens for a broad class of initial data - specifically, near the points where $\omega_0(x)$ vanishes and
the real part of $H \omega_0$ has the right sign.

A more general model has been proposed by De Gregorio \cite{DG1}, \cite{DG2}:
\begin{equation}\label{111}
\partial_t \omega + u \partial_x \omega = \omega H \omega,\,\,\,u_x = H\omega, \,\,\,\omega(x,0) = \omega_0(x).
\end{equation}
This model includes the advection term. Amazingly, the question of whether the solutions to \eqref{111} are globally
regular or can blow up in finite time is currently open. Numerical simulations appear to suggest global regularity \cite{OSW}, but the
mechanism for it is not well understood. Recently, global regularity near a manifold of equilibria as well as other interesting
features of the solutions of \eqref{111} have been shown in \cite{JSS}. Variants of \eqref{111} and other related models appear in
for example \cite{BKP,CasCor,EJ3,EK1,Wunsch}, where further references can be found.

Already in \cite{HouLuo}, Hou and Luo proposed a simplified one-dimensional model specifically designed
to gain insight into the singularity formation process in the scenario described in Section~\ref{hyp}.
This model is given by
\begin{eqnarray}
\partial_t \omega + u \partial_x \omega = \partial_x \theta, \nonumber \\
\partial_t \theta + u \partial_x \theta = 0, \,\,\, u_x = H \omega.
\label{1dhl}
\end{eqnarray}
Here as above $H$ is the Hilbert transform, and the setting can be either periodic or the entire axis with some decay of the initial data.
The model \eqref{1dhl} can be thought of as an effective equation on the $x_2=0$ axis in the Boussinesq case
(see \eqref{boussinesq} and Figure~\ref{Bousfig}) or on the boundary of the cylinder in the 3D axi-symmetric Euler case.
The model can be derived from the full equations under certain boundary layer assumption:
that $\omega(x,t)$ is concentrated in a boundary layer of width $a$ near $x_2=0$ axis and is independent of $x_2$,
that is
$\omega(x_1,x_2,t) = \omega(x_1,t) \chi_{[0,a]}(x_2).$
Such assumption is necessary to close the equation and reduces the half-plane Biot-Savart law to $u_x=H\omega$ in the main order;
the parameter $a$ enters into the additional term that is non-singular and is dropped from \eqref{1dhl}.
See \cite{HouLuo}, \cite{CHKLSY} for more details. We will call the system \eqref{1dhl} the HL model.

The HL model is still fully nonlocal. A further simplification was proposed in \cite{CKY}, where the Biot-Savart law
has been replaced with
\be\label{ckybs}
u(x,t) = - x \int_x^1 \frac{\omega(y,t)}{y}\,dy.
\ee
Here the most natural setting is on an interval $[0,1]$ with smooth initial data supported away from the endpoints.
The law \eqref{ckybs} is motivated by the velocity representation in Lemma~\ref{mainlemma} above, as it is the simplest one dimensional
analog of such representation. This law is ``almost local" - if one divides $u$ by $x$ and differentiates, one gets local expression.
We will call the model \eqref{ckybs} the CKY model.

For both HL and CKY models, local well-posedness in a reasonable family of spaces (such as sufficiently regular Sobolev spaces) is not
difficult to obtain. In \cite{CKY}, finite time blow up has been proved for the CKY model. The proof used analysis of the trajectories
and of the nonlinear feedback loop generated by the forcing term $\partial_x \theta$. We will sketch a very similar argument below. The proof does not provide
a detailed blow picture. In a later work \cite{LH}, more precise picture of blow up was established with aid of computer
assisted proof. It shows self-similar  behavior near the origin properly matched with the outside region. For the original HL model, finite
time blow up has been proved in \cite{CHKLSY}. For the model including the additional term obtained from the boundary layer assumption
into Biot-Savart law, finite time blow up has been proved in \cite{DKX}. We now sketch a variant of the blow up proof for the HL model \eqref{1dhl}.

Let us consider an HL model on $[0,L]$ with periodic boundary conditions. In this setting, using the expression for periodic Hilbert transform,
the Biot-Savart law becomes
\[ u_x(x) = H \omega(x) = \frac1L P.V. \int_0^L \omega(y) \cot[\mu (x-y)]\,dy, \]
where $\mu = \pi/L.$ Integration leads to
\be\lb{perBS}
u(x) = \frac{1}{\pi} \int_0^L \omega(y) \log|\sin[\mu(x-y)]|\,dy.
\ee
The initial data will be chosen as follows: $\omega_0$ is odd, which together with periodicity implies that it is also
odd with respect to $x=L/2,$ and satisfy $\omega_0(x) \geq 0$ if $x \in [0,L/2].$ The initial density $\theta_0$ is even
with respect to both $0$ and $L/2,$ and satisfies $\theta_0' \geq 0$ for $x \in [0,L/2].$ The solution, while it exists, will satisfy
the same properties. The symmetry assumptions on $\omega$ lead to the following version of the Biot-Savart law, which can be verified
by direct computation.
\begin{lemma}\label{l3a}
Let $\omega$ be periodic with period $L$ and odd at $x = 0$ and let $u$ be defined by \eqref{perBS}. Then for any $x \in [0,\frac{1}{2} L]$,
\begin{equation}
  u(x) \cot(\mu x) = -\frac{1}{\pi} \int_{0}^{L/2} K(x,y) \omega(y) \cot(\mu y)\,dy,
  \label{eqn_u_ker}
\end{equation}
where
\begin{equation}
  K(x,y) = s \log \biggl| \frac{s+1}{s-1} \biggr|\qquad \text{with}\qquad s = s(x,y) = \frac{\tan(\mu y)}{\tan(\mu x)}.
  \label{eqn_K}
\end{equation}
Furthermore, the kernel $K(x,y)$ has the following properties:
\begin{enumerate}
  \item $K(x,y) \geq 0$ for all $x,\, y \in (0,\frac{1}{2} L)$ with $x \neq y$;
  \vspace{1mm}
  \item $K(x,y) \geq 2$ and $K_{x}(x,y) \geq 0$ for all $0 < x < y < \frac{1}{2} L$;
\end{enumerate}
\label{lmm_kernel}
\end{lemma}
The key observation is a certain positivity property that will help us control the behavior of trajectories \cite{CHKLSY}.
\begin{lemma}\label{l4}
Let the assumptions in Lemma \ref{lmm_kernel} be satisfied and assume in addition that $\omega \geq 0$ on $[0,\frac{1}{2} L]$. Then for any $a \in [0,\frac{1}{2} L]$,
\begin{equation}
  \int_{a}^{L/2} \omega(x) \bigl[ u(x) \cot(\mu x) \bigr]_{x}\,dx \geq 0.
  \label{eqn_wvy_pos}
\end{equation}
\label{lmm_wvy_pos}
\end{lemma}
With these two lemmas, the rest of the proof proceeds as follows.
Towards a contradiction, let us assume that there exists a global solution $(\theta,\omega)$ to \eqref{1dhl} with
the initial data as described in the beginning of this section and
 denote $A := \theta_{0}(\frac{1}{2} L) > 0$. Since $\theta_{0}$ is assumed to be increasing on $[0,\frac{1}{2} L]$, we can choose a decreasing sequence $x_{n}$ in $(0,\frac{1}{2} L)$ with
$n \geq 0$ such that $\theta_{0}(x_{n}) = [2^{-1} + 2^{-(n+2)}] A$. Note that $x_{0} < \frac{1}{2} L$ since $\theta_{0}(x_{0}) < \theta_{0}(\frac{1}{2} L)$.

For $x_{n}$ defined as above, let $\Phi_{n}(t)$ denote the characteristics of \eqref{1dhl} originating from $x_{n}$, that is, let $\frac{d}{dt} \Phi_{n}(t) = u(\Phi_{n}(t),t)$ with $\Phi_{n}(0) = x_{n}$.
Lemma \ref{lmm_kernel} then implies the following estimate on the evolution of $\Phi_{n}$:
\begin{equation}
  \frac{d}{dt} \Phi_{n}(t) = u(\Phi_{n}(t),t) \leq -\frac{2}{\pi} \tan(\mu \Phi_{n}(t)) \int_{\Phi_{n}(t)}^{L/2} \omega(y,t) \cot(\mu y)\,dy \leq -\frac{2\mu}{\pi}\, \Phi_{n}(t) \Omega_{n}(t),
  \label{eqn_dPhi}
\end{equation}
where for simplicity we have written
\begin{equation}
  \Omega_{n}(t) := \int_{\Phi_{n}(t)}^{L/2} \omega(y,t) \cot(\mu y)\,dy.
  \label{eqn_def_Omege}
\end{equation}
Introducing the new variable $\psi_{n}(t) := -\log \Phi_{n}(t)$, we may write \eqref{eqn_dPhi} as
\begin{equation}
  \frac{d}{dt} \psi_{n}(t) \geq \frac{2\mu}{\pi}\, \Omega_{n}(t).
  \label{eqn_dpsi}
\end{equation}
Then for each $n \geq 1,$ we have
\begin{align*}
  \frac{d}{dt} \Omega_{n}(t) & = \int_{\Phi_{n}(t)}^{L/2} \omega(y,t) \bigl[ u(y,t) \cot(\mu y) \bigr]_{y}\,dy + \int_{\Phi_{n}(t)}^{L/2} \theta_{y}(y,t) \cot(\mu y)\,dy \\
  & \geq \int_{\Phi_{n}(t)}^{\Phi_{n-1}(t)} \theta_{y}(y,t) \cot(\mu y)\,dy \\
  & \geq \cot(\mu \Phi_{n-1}(t)) \bigl[ \theta_{0}(x_{n-1}) - \theta_{0}(x_{n}) \bigr] = 2^{-(n+2)} A \cot(\mu \Phi_{n-1}(t)),
\end{align*}
where in the second step we have used Lemma \ref{lmm_wvy_pos} and the fact that $\theta_{x} \geq 0$ on $[0,\frac{1}{2} L]$. To find a lower bound for the right hand side, note that for any
fixed $z \in (0,\frac{1}{2} \pi)$, there exists some constant $c > 0$ depending only on $z$ such that $\cot(x) > cx^{-1}$ for any $x \in (0,z]$. In our situation, we have
$\mu \Phi_{n-1}(t) \leq \mu \Phi_{0}(t) \leq \mu x_{0} < \frac{1}{2} \pi$, and as a result there exists some constant $c_{0} > 0$ depending only on $\mu$ and $x_{0}$ such that
$\cot(\mu \Phi_{n-1}(t)) \geq c_{0} [\Phi_{n-1}(t)]^{-1}$. This leads to the estimate
\begin{equation}
  \frac{d}{dt} \Omega_{n}(t) \geq 2^{-(n+2)} c_{0} A e^{\psi_{n-1}(t)}.
  \label{eqn_dOmega}
\end{equation}
Once we have \eqref{eqn_dpsi}, \eqref{eqn_dOmega}, the proof of finite time blow is fairly straightforward. Details of a similar argument can be found in \cite{CKY}.
We can choose $A$ large enough to show inductively that
$\psi_n(t_n) \geq b n + a$ for some suitably chosen $a \in \Rm,$ $b >0$ and an increasing sequence $t_n \rightarrow T < \infty.$ This implies that $\theta$ has to develop a shock at
$x=0$ by time $T$ - unless blow up happens before that in some other fashion (invalidating regularity assumptions underlying our estimates such as integration by parts).
In fact, the informal flavor of bounds \eqref{eqn_dpsi}, \eqref{eqn_dOmega} is that of $F'' \geq ce^{cF}$ differential inequality, leading to dramatically fast growth.

\section{The SQG patch problem: a blow up blueprint}\label{patch}

Part of the difficulty in securing rigorous understanding of the Hou-Luo blow up scenario for 3D axi-symmetric Euler or 2D inviscid Boussinesq
equation lies in growth of $\omega,$ which destroys the estimate of error terms in Lemma~\ref{mainlemma}. The error terms may no longer be
of smaller order than the main term. In fact, heuristic computations taking $\omega$ that behaves approximately like some inverse power of $x_1$ in a certain
region near origin - an ansatz that appears to be in agreement with the numerical simulations - indicate that the error terms will now be
of the same order as the main term obtained by integration over the bulk. In this section, we discuss a different setting in which this
situation is the case - the portion of the Biot-Savart integral pushing the solution towards blow up has the same order as the part pushing
in the opposing direction. Nevertheless, the conclusion is finite time blow, essentially due to presence of a parameter that can be
used to overcome the error term. The setting is that of modified SQG patch solutions in the half-plane.

Namely, let us in this section set $D=\Rm^2_+=\{(x_1,x_2) \left| \,x_2 \geq 0 \right. \}.$
The Bio-Savart law for the patch evolution on the half-plane $D := \Rm\times\Rm^+$~is
\[
u=\nabla^\perp (-\Delta_D)^{-1+\alpha}\omega,
\]
with $\Delta_D$ being the Dirichlet Laplacian on $D$, which can also be written as
\begin{equation}
u(x, t) :=  \int_D \left( \frac{(x-y)^\perp}{|x-y|^{2+2\alpha}} -
\frac{(x-\bar y)^\perp}{|x-\bar y|^{2+2\alpha}} \right) \omega(y,t) dy
\label{eq:velocity_law}
\end{equation}
The case $\alpha = 0$ corresponds to the 2D Euler equation, while $\alpha=1/2$ to the SQG equation; the range $0 \leq \alpha \leq 1$
is called modified SQG. Note that $u$ is divergence free and tangential to the boundary.
A traditional approach to the 2D Euler ($\alpha=0$) vortex patch evolution, going back to Yudovich (see~\cite{MP}
for an exposition) is via the corresponding  flow map. The active scalar $\omega$ is advected by~$u$ from (\ref{eq:velocity_law}) via
\begin{equation}\label{1.31}
\omega(x,t) = \omega \left(\Phi^{-1}_t(x),0\right),
\end{equation}
where
\begin{equation}\label{eq:alpha}
\frac{d}{dt}\Phi_t(x) = u\left(\Phi_t(x),t \right) \qquad\text{and}\qquad  \Phi_0(x)=x.
\end{equation}
The initial condition $\omega_0$ for  \eqref{eq:velocity_law}-\eqref{eq:alpha} is patch-like,
\begin{equation}\label{patchlikeid}
\omega_0 = \sum_{k=1}^N \theta_k \chi_{\Omega_{0k}},
\end{equation}
with $\theta_1,\dots,\theta_N\neq 0$ and $\Omega_{01},\dots,\Omega_{0N}\subseteq D$  bounded open sets, whose closures $\overline{\Omega_{0,k}}$ are
pairwise disjoint and whose boundaries~$\partial\Omega_{0k}$
are simple closed curves. The question of regularity of solution in patch setting becomes the question of the conservation of regularity class
of the patch boundary, as well as lack of self-intersection or collisions between different patches.

One reason the Yudovich theory works for the 2D Euler equations is that for $\omega$ which is (uniformly in time) in $L^1\cap L^\infty$, the velocity field
$u$ given by (\ref{eq:velocity_law}) with $\alpha=0$ is log-Lipschitz in space, and the flow map $\Phi_t$ is  everywhere well-defined (see e.g. \cite{MB,MP}).
In our situation,
when $\omega$ is a patch solution and $\alpha>0$, the flow $u$ from (\ref{eq:velocity_law}) is smooth away from the patch boundaries  $\partial\Omega_k(t)$ but is only
$1-\alpha$ H\"older constinuous at $\partial\Omega_k(t),$ which is exactly where one needs to use the flow map.
This creates significant technical difficulties in proving
local well-posedness of patch evolution in some reasonable functional space. For the case without boundaries, local well-posedness has been proved in \cite{Rodrigo}
for $C^\infty$ patches and for Sobolev $H^3$ patches in \cite{g} for $0<\alpha \leq 1.$ A naive intuition on why patch evolution can be locally well-posed for $\alpha >0$ is that
the below-Lipschitz loss of regularity only affects the tangential component of the fluid velocity at patch boundary. The normal to patch component, that intuitively should determine
the evolution of the patch, retains stronger regularity.

In presence of boundaries, the problem is harder. Intuitively, one reason for the difficulties can be explained as follows. In the simplest case of half-plane the reflection principle
implies that the boundary can be replaced by a reflected patch (or patches) of the opposite sign. If the patch is touching the boundary, then the reflected and original
patch are touching each other, and the low regularity tangential component of the velocity field generated by the reflected patch has strong influence on the boundary of the original patch
near touch points.
Even in the 2D Euler case, the global regularity for patches in general domains with boundaries is currently open (partial results for patches not touching the boundary
or with loss of regularity can be found in \cite{d}, \cite{d2}).
In the half-plane, a global regularity result has been recently established in \cite{KRYZ}:
\begin{theorem}\label{thmeuler1}
Let $\alpha =0$ and $\gamma\in(0,1]$.  Then for each $C^{1,\gamma}$ patch-like initial data
$\omega_0$,  there exists a unique global $C^{1,\gamma}$ patch solution $\omega$ to \eqref{1.31}, \eqref{eq:velocity_law}, \eqref{eq:alpha}
with $\omega(\cdot,0)=\omega_0$.
\end{theorem}
In the case $\frac{1}{24}> \alpha>0$ with boundary, even local well-posedness results are highly non-trivial. The following result has been proved in \cite{KYZ}
for the half-plane.
 \begin{theorem}\label{T.1.1}
If  $\alpha\in(0,\frac 1{24})$, then for each $H^3$ patch-like initial data $\omega_0$,  there exists a unique local $H^3$
patch solution $\omega$ with $\omega(\cdot,0)=\omega_0$.  Moreover, if the maximal time
$T_\omega$ of existence of $\omega$ is finite, then at $T_\omega$ a singularity forms: either two patches touch, or a patch boundary touches
itself or  loses $H^3$ regularity.
\end{theorem}
We note that one has to be careful in the definition of solutions in this case as trajectories \eqref{eq:alpha} may not be unique.
Solutions can be defined in a weak sense by pairing with a test function, or in an appealing geometric way by specifying evolution
of patch boundary with velocity \eqref{eq:velocity_law} in the sense of Hausdorff distance; see \cite{KYZ} for more details.
The constraint $\alpha < \frac{1}{24}$ appears due to estimates near boundary; it is not clear if it is sharp.


On the other hand, in \cite{KRYZ}, it was proved that for any $\alpha>0,$ there exist patch-like initial data leading to finite time blow up.
\begin{theorem}\label{main1234}
Let $\alpha\in(0,\frac 1{24})$.  Then there are $H^3$ patch-like initial data $\omega_0$ for which the unique local $H^3$
patch solution $\omega$ with $\omega(\cdot,0)=\omega_0$
becomes singular in finite time (i.e., its maximal time of existence $T_\omega$ is finite).
\end{theorem}
Together, Theorems~\ref{thmeuler1} and \ref{main1234} give rigorous meaning to calling the 2D Euler equation critical. In the half-plane patch framework
$\alpha=0$ is the exact threshold for phase transition from global regularity to possibility of finite time blow up.

In what follows, we will sketch proof the blow up Theorem~\ref{main1234}.
We concentrate on the main ideas only; full details can be found in \cite{KRYZ}.
Let us describe the initial data set up.
Denote $\Omega_1:=(\eps,4)\times(0,4)$, $\Omega_2:=(2\eps,3)\times(0,3)$, and
let $\Omega_0 \subseteq D^+ \equiv \Rm^+ \times \Rm^+$ be an open set whose boundary is a smooth simple closed curve and
which satisfies~$\Omega_2 \subseteq \Omega_0 \subseteq \Omega_1$. Here $\epsilon$ is a small parameter
depending on $\alpha$ that will be chosen later.

Let $\omega(x,t)$ be the  unique $H^3$ patch solution corresponding to the initial data
\begin{equation}
\omega(x,0) := \chi_{\Omega_0}(x) - \chi_{\tilde{\Omega}_0}(x)
\label{def:omega_0}
\end{equation}
with maximal time of existence $T_\omega>0$.
Here, $\tilde{\Omega}_0$ is the reflection of $\Omega_0$ with respect to the $x_2$-axis.
Then
\begin{equation} \label{6.0}
\omega(x, t) = \chi_{\Omega(t)}(x) - \chi_{\tilde{\Omega}(t)}(x)
\end{equation}
for $t\in[0,T_\omega)$, with $\Omega (t):=\Phi_t(\Omega_0).$ 
It can be seen from \eqref{eq:velocity_law} that the rightmost point of the left patch on the $x_1$-axis and
the leftmost point of the right patch on the $x_1$-axis will move toward each other. 
In the case of the 2D Euler equations $\alpha=0$, Theorem \ref{thmeuler1} shows that the two points never reach the origin.
When $\alpha>0$ is small, however, it is possible to control the evolution sufficiently well to show that --- unless the solution develops another
singularity earlier --- both points will reach the origin in a finite time. The argument yielding such control is fairly subtle, and
the estimates do not extend to all $\alpha<\frac 12$, even though one would expect singularity formation to persist for more singular equations.
This situation is not uncommon in the field: there is plenty of examples with the infinite in time growth of derivatives for the
smooth solutions of 2D Euler equation, while none are available for the more singular SQG equation \cite{KN5}.

\begin{figure}
\begin{center}
\includegraphics[scale=1]{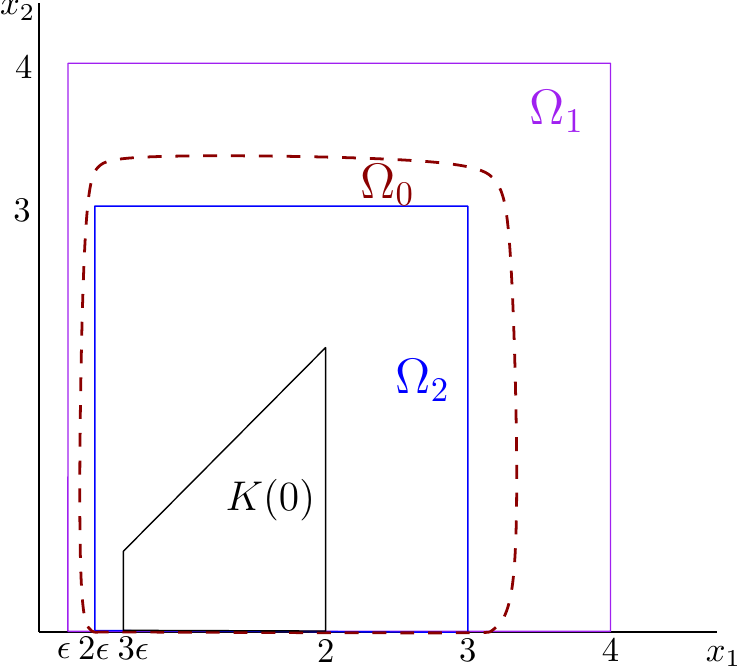}
\end{center}
\caption{The domains $\Omega_1, \Omega_2, \Omega_0$, and $K(0)$ (with $\omega_0=\chi_{\Omega_0}-\chi_{\tilde  \Omega_0}$). \label{fig:def_domains}}
\end{figure}

To show finite time blow up, we will deploy a barrier argument. Define
\begin{equation}\label {6.2}
K(t):=\{ x\in D^+ \,:\, x_1\in(X(t),2) \text{ and } x_2\in(0,x_1)\}
\end{equation}
for $t\in[0,T]$, with $X(0) =3\epsilon.$ Clearly, $K(0) \subset \Omega(0).$
Set the evolution of the barrier by
\begin{equation}\label{aug2110}
X'(t)=-\frac 1{100\alpha} X(t)^{1-2\alpha}.
\end{equation}
Then $X(T)=0$ for $T = 50(3\epsilon)^{2\alpha}.$ So if we can show that $K(t)$ stays inside $\Omega(t)$
while the patch solution stays regular, then we obtain that singularity must form by time $T:$
the different patch components will touch at the origin by this time unless regularity is lost before that.

The key step in the proof involves estimates of the velocity near origin.
In particular, $u_1$ needs to be sufficiently negative to exceed
the barrier speed \eqref{aug2110}; $u_2$ needs to be sufficiently positive in order to ensure
that $\Omega (t)$ cannot cross the barrier along its diagonal part. Note that it suffices
to consider the part of the barrier that is very close to the origin, on the order $\sim \epsilon^{2\alpha}.$
Indeed, the time $T$ of barrier arrival at the origin has this order, and the fluid velocity
satisfies uniform $L^\infty$ bound that follows by a simple estimate which uses only $\alpha < 1/2.$
Thus the patch $\Omega(t)$ has no time to reach more distant boundary points of the barrier before
formation of singularity.

Let us focus on the estimates for $u_1.$
For $y=(y_1,y_2)\in \bar D^+ = \Rm^+\times\Rm^+$, we denote $\bar y:=(y_1,-y_2)$
and $\tilde y:=(-y_1,y_2).$
Due to odd symmetry, \eqref{eq:velocity_law} becomes (we drop $t$ from the notation in this sub-section)
\begin{equation}\label{BS12}
u_1(x) = -\int_{D^+} K_1(x,y) \omega(y) dy,
\end{equation}
where
\begin{equation}\label{BS1det}
K_1(x,y) =
\underbrace{\frac{y_2-x_2}{|x-y|^{2+2\alpha}}}_{K_{11}(x,y)} -
\underbrace{\frac{y_2-x_2}{|x-\tilde y|^{2+2\alpha}}}_{K_{12}(x,y)} -
 \underbrace{\frac{y_2+x_2}{|x+y|^{2+2\alpha}}}_{K_{13}(x,y)} +
\underbrace{\frac{y_2+x_2}{|x-\bar y|^{2+2\alpha}}}_{K_{14}(x,y)},
\end{equation}

Analyzing \eqref{BS1det}, it is not hard to see that we can split the region of integration in the Biot-Savart
law according to whether it helps or opposes the bounds we seek. Define
\[
u_{1}^{bad}(x) := - \int_{\Rm^+\times(0,x_2)} K_1(x,y) \omega(y) dy \,\,\,
{\rm and} \,\,\,
u_{1}^{good}(x) := - \int_{\Rm^+\times(x_2,\infty)} K_1(x,y) \omega(y) dy.
\]

The following two lemmas contain key estimates.

\begin{lemma}[Bad part]
Let $\alpha\in(0,\frac 12)$ and assume that $\omega$ is odd in $x_1$ and $0\le \omega\le 1$ on $D^+$.
If $x\in \overline{D^+}$ and $x_2\leq x_1$, then
\be\lb{bades}
u_1^{bad}(x) \leq \frac{1}{\alpha}\left(\frac{1}{1-2\alpha} -2^{-\alpha}\right) x_1^{1-2\alpha}.
\ee
\label{lemma:bad_parts_u1}
\end{lemma}
The proof of this lemma uses \eqref{BS1det} and after cancellations leads to the bound
\be\lb{bades1} u_1^{bad}(x) \leq -\int_{(0,2x_1)\times(0,x_2)} \frac{y_2-x_2}{|x-y|^{2+2\alpha}}  dy, \ee
which gives \eqref{bades}

In the estimate of the good part, we need to use a lower bound on $\omega$ that will be provided by the barrier. Define
\begin{equation}
A(x) := \left\{ y\,:\,  y_1\in \left(x_1,x_1+1 \right) \text{ and } y_2 \in(x_2, x_2 + y_1 - x_1) \right\}.
\label{eq:def_A}
\end{equation}
\begin{lemma}[Good part]
\label{lemma:u1_good}
Let $\alpha\in(0,\frac 12)$ and assume that $\omega$ is odd in $x_1$ and for some $x\in  \overline{D^+}$
we have $\omega\ge \chi_{A(x)}$ on $D^+$, with $A(x)$ from \eqref{eq:def_A}.
There exists $\delta_\alpha\in(0,1)$, depending only on $\alpha$, such that the following holds.
\\[0.2cm]
If $x_1\le \delta_\alpha$, then
\begin{equation*}
u_1^{good}(x) \leq -\frac{1}{6\cdot 20^\alpha \alpha}x_1^{1-2\alpha}.
\end{equation*}
\end{lemma}

Here analysis of \eqref{BS1det} leads to
\[
u_1^{good}(x) \leq -\underbrace{\int_{A_1}\frac{y_2-x_2}{|x-y|^{2+2\alpha}}  dy}_{T_1}
+ \underbrace{ \int_{A_2}\frac{y_2-x_2}{|x-y|^{2+2\alpha}}  dy}_{T_2},
\]
with the domains
\begin{equation*}
\begin{split}
A_1 &:= \left\{ y \,:\, y_2 \in \left(x_2, x_2 + 1 \right) \text{ and } y_1\in \left(x_1+y_2-x_2, 3x_1+y_2-x_2 \right) \right\},\\
A_2 &:=  \left(x_1+1, 3x_1+1\right) \times \left(x_2, x_2 + 1 \right).
\end{split}
\end{equation*}
The term $T_2$ can be estimated by $Cx_1,$ since the region of integration $A_2$ lies at a distance $\sim 1$ from the singularity.
A relatively direct estimate of the term $T_1$ leads to the result of the Lemma.

A distinctive feature of the problem is that estimates for the ``bad'' and ``good'' terms appearing in Lemmas~\ref{lemma:bad_parts_u1} and \ref{lemma:u1_good}
above have the same order of magnitude $x_1^{1-2\alpha}.$ This is unlike the 2D Euler double exponential growth construction, where we were able to isolate the main term.
To understand the balance in the estimates for the ``bad'' and ``good'' terms, note that the ``bad'' term estimate comes from integration
of the Biot-Savart kernel over rectangle $(0,2x_1) \times (0,x_2),$ while the good term estimate from integration of the same kernel over the region
$A_1$ above. When $\alpha$ is close to zero, the kernel is longer range, and the more extended nature of the region $A_1$ makes the ``good'' term dominate.
In particular, the coefficient  $\frac{1}{\alpha}\left(\frac{1}{1-2\alpha}-2^{-\alpha}\right)$ in front of $x_1^{1-2\alpha}$ in Lemma~\ref{lemma:bad_parts_u1} converges to
to finite limit as $\alpha \rightarrow 0,$ while the coefficient $\frac{1}{6\cdot 20^\alpha \alpha}$ in Lemma~\ref{lemma:u1_good} tends to infinity. On the other hand,
when $\alpha \rightarrow \frac12,$
the singularity in the Biot-Savart kernel is strong and getting close to non-integrable. Then it becomes important that the ``bad'' term integration region
contains an angle $\pi$ range near the singularity, while the ``good'' region only $\frac{\pi}{4}.$ For this reason, controlling the ``bad'' term for larger values of $\alpha$ is problematic - although there is no reason why there cannot be a different, more clever argument achieving this goal.

It is straightforward to check that the dominance of the ``good'' term over ``bad'' one extends to the range $\alpha \in (0,\frac{1}{24}),$
and that in this range we get as a result
\be\label{u1dec}
u_1(x,t) \leq -\frac{1}{50\alpha} x_1^{1-2\alpha}
\ee
for $x=(x_1,x_2)$ such that $x_1 \leq \delta_\alpha$ and $x_1 \geq x_2.$ A similar bound can be proved showing that
\be\label{u2dec}
u_2(x,t) \geq \frac{1}{50\alpha} x_2^{1-2\alpha}
\ee
for $x=(x_1,x_2)$ such that $x_2 \leq \delta_\alpha$ and $x_1 \leq x_2.$

\begin{figure}[htbp]
\begin{center}
\includegraphics[scale=1.2]{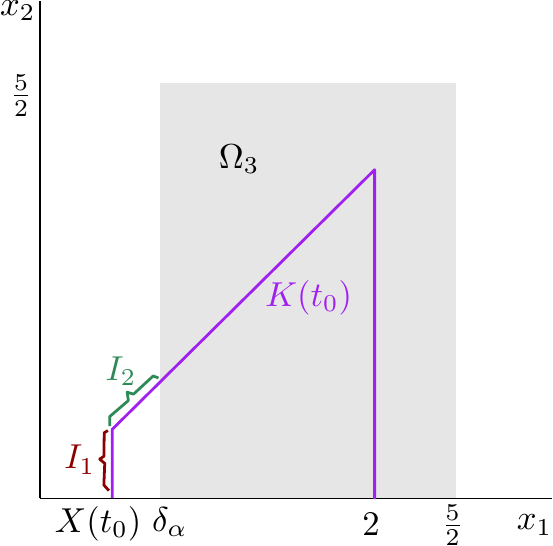}
\caption{The segments $I_1$ and $I_2$ and the sets $\Omega_3$ and $K(t_0)$.\label{fig:t1}}
\end{center}
\end{figure}

The proof is completed by a contradiction argument, where we assume that the barrier $K(t)$ catches up with the patch $\Omega(t)$
at some time $t=t_0 < T$ of first contact. Taking $\epsilon$ sufficiently small compared to $\delta_\alpha$ from Lemma~\ref{lemma:u1_good}, we can make sure
the contact can only happen on the intervals $I_1$ and $I_2$ along the boundary of the barrier $K(t_0)$ appearing on Figure~\ref{fig:t1}.
But then bounds \eqref{u1dec}, \eqref{u2dec} and the evolution of the barrier prescription \eqref{aug2110} lead to the conclusion that
the barrier should have been crossed at $t < t_0,$ yielding a contradiction; full details can be
found in \cite{KRYZ}).


\section{Discussion}\label{disc}

There are a few more recent papers that have contributed towards understanding the hyperbolic point blow up scenario. Two-dimensional simplified models
of the 2D Boussinesq system have been considered in \cite{HORY} and in \cite{KT}. In both cases, the derivative forcing term in \eqref{boussinesq}
is replaced by a simpler sign-definite approximation $\frac{\theta}{x_1},$ and the Biot-Savart law
is replaced by a simpler version $u = (-x_1\Omega(x,t),x_2 \Omega(x,t)).$ In \cite{HORY}, $\Omega$ takes form similar to the 2D Euler example
\eqref{omegamain}. In \cite{KT}, $\Omega$ is closely related but is also chosen to keep $u$ incompressible. Both papers prove finite time blow up,
\cite{HORY} by a sort of barrier argument while the argument \cite{KT} deploys an appropriate Lyapunov-type functional.

A very interesting recent work by Elgindi and Jeong takes a different approach \cite{EJ1}, \cite{EJ2}. In \cite{EJ1}, they look at a class of scale invariant
solutions for the 2D Boussinesq system that satisfy $\frac{1}{\lambda} u(\lambda x,t) = u(x,t)$ and
$\frac{1}{\lambda} \theta(\lambda x,t) = \theta(x,t).$ Observe that this class allows velocity and density that grow linearly at infinity.
Also, the solution is not regular at the origin: for example the vorticity is just $L^\infty.$
The setting is a sector which has size $\frac{\pi}{2}$ (and some results can be generalized to other angles $<\pi$).
First, they prove a local well-posedness theorem in a class of solutions that includes scale invariant solutions; additional symmetry assumptions
are needed for this result. Secondly, for such solutions, they obtain an effective one-dimensional equation, some solutions of which
are shown to blow up in finite time. These are not the first examples of infinite energy solutions (see \cite{Childr}, \cite{Const},
and \cite{jWu}). However in \cite{EJ1} a procedure to cut off the solution at infinity while maintaining finite time blow up property
is carried out. This yields finite energy solutions leading to finite time blow up - in the sense that $\int_0^T
\|\nabla u(\cdot, t)\|_{L^\infty}\,dt \rightarrow \infty $ at blow up time. The vertex of the sector is a hyperbolic stagnation point
of the flow, making connection to the Hou-Luo scenario. 
In \cite{EJ2}, similar results are proved in the 3D axi-symmetric
Euler case; here the domain is given by $(1+\epsilon |z|)^2 \leq |x|^2+|y|^2$ with an arbitrary $\epsilon >0.$ 

The main challenge to analyzing smooth solutions to 2D Boussinesq and 3D axi-symmetric Euler equations in this context stems from difficulties estimating
the velocity produced by Biot-Savart law with growing vorticity. It does not appear that there is a clear main term, as in 2D Euler example,
or a clear small parameter to play in the same order of magnitude opposing terms, as in modified SQG patches. On the other hand, all model problems
point to finite time blow up outcome in the original Hou-Luo scenario. The challenge is finding enough controllable structures to
carry through rigorous analysis.

\end{document}